\documentclass[10pt]{amsart}
\usepackage{amsfonts} 
\usepackage{amscd}
\usepackage{amssymb}
\usepackage{float}
\usepackage{amsmath}
\usepackage{graphicx}
\usepackage{amsxtra}
\usepackage{amsthm,latexsym,amstext,delarray,
amsmath,epsfig}
\usepackage{eucal}
\usepackage[all]{xy}

\usepackage[OT2,OT1]{fontenc}
\DeclareMathAlphabet{\cyr}{OT2}{wncyr}{m}{n}

\newcommand{\ring}[1]{\mathbb{#1}}
\newcommand{\C}{\ring{C}} \newcommand{\Q}{\ring{Q}}
\newcommand{\Pe}{\ring{P}}
\newcommand{\Z}{\ring{Z}} 
 \newcommand{\G}{\ring{G}} 
\newcommand{\F}{\ring{F}} 
\newcommand{\bu}{\bullet}

\newcommand{\ga}{\ring{G}_a}
\newcommand{\gm}{\ring{G}_m}
\newcommand{\be}{\begin{equation}}
\newcommand{\ee}{\end{equation}}
\newcommand{\nd}{\noindent}
\newcommand{\End}{\operatorname{End}}
\newcommand{\Ext}{\operatorname{Ext}}

\newcommand{\Coker}{\operatorname{Coker}}
\newcommand{\Ker}{\operatorname{Ker}}
\def\Hom{\textup{Hom}}
\def\1{{\mu\mkern-6mu\mu}}
\newcommand{\Sha}{\cyr{X}}
\def\half{{\frac{1}{2}}}

\begin{document}

\title[Values at \(s=\half\)]{Values of zeta functions at \(s=\half\)}
\author[Ramachandran]{Niranjan Ramachandran}
\dedicatory{To PVR and SS, in gratitude}
\thanks{Partially supported by MPIM, IHES, and a GRB Summer Fellowship
  from the University of Maryland}
\address{N.~Ramachandran: Department of Mathematics \\ University of
Maryland \\ College Park, MD 20912 USA}
\address{Institut des Hautes \'{E}tudes Scientifiques, Route de
  Chartres, F-91440, Bures-sur-Yvette, France}
\email{atma@math.umd.edu}
\date{\today}

\theoremstyle{plain}
\newtheorem{thm}{Theorem}
\newtheorem{lem}[thm]{Lemma}
\newtheorem{cor}[thm]{Corollary}
\newtheorem{prop}[thm]{Proposition}
\newtheorem{conj}[thm]{Conjecture}
\newtheorem{quest}[thm]{Question}
 
\theoremstyle{definition}
\newtheorem{rem}[thm]{Remark}
\newtheorem{defn}[thm]{Definition}
\newtheorem{ex}[thm]{Example}
\newtheorem{obser}[thm]{Observation}

\begin{abstract} 
We  study the
behaviour near \(s=\half\) of zeta functions of varieties over finite
fields \(\F_q\)  with \(q\) a square. The main result is an
Euler-characteristic formula for the
square of the special value 
at \(s=\half\).  The Euler-characteristic is constructed from the
Weil-\'etale cohomology of a certain supersingular elliptic curve.  
\end{abstract}   
\maketitle

\section*{Introduction}

Let \(V\) be a integral scheme of finite type over Spec \(\Z\). The special
values \(a_V(n)\) at integers \(s=n\) of the zeta
function \(\zeta(V,s)\) of \(V\)  are conjecturally related to deep arithmetical
invariants of \(V\). One may ask if the special values of
\(\zeta(V,s)\) at non-integral values of \(s\), e.g. \(s=\half\),  also admit an
arithmetical interpretation.   The simplest example -- and perhaps
the most interesting -- arises from 
the central value of the Riemann zeta function \(\zeta(s) = 
\zeta({\rm  Spec}~ \Z,s)\): 
\begin{center}{\em Is there a motivic interpretation of 
  \(\zeta(\half)\)?}\end{center}
 
For instance, it is unknown if \(\zeta(\half)\) is a period in the
sense of \cite{kz}; however, for certain triple product L-functions,
the work of M.~ Harris and S.~Kudla 
\cite{hk} relates the central critical value (at \(s=\half\)) to
periods, cf. \cite[Eq.~(48), p.459]{sa}. 

The motivic 
philosophy indicates that \(a_V(n)\) depends on the interaction of the
motive \(h(V)\) of \(V\) with \(\Z(\pm n)\) (power of the Tate
motive). This leads us to suspect that the
special value  at \(s=\half\) is governed by an
exotic object (unknown to exist, as yet): a square root \(\Z(\half)\)
of the Tate motive \(\Z(1)\) over Spec~\(\Z\). 

It is clear that the investigation of special values at \(s=\half\) should begin with
the important case of varieties over finite fields. Namely, consider
the zeta function \(\zeta(X,s)\) of a smooth projective variety \(X\)
over a finite field 
\(k= \F_q\) of characteristic \(p >0\).  The function \(Z(X,t)\), defined by \(Z(X,q^{-s}) :=
\zeta(X,s)\),  is a rational function of \(t\) with integer
coefficients.  For any integer \(n\ge 0\), the order of the pole 
\(\rho_n:= -{\rm ord}_{s=n} \zeta(X,s)\) at \(s=n\) and the  special value
\(a_X(n)\) of \(Z(X,t)\) at  \(t=q^{-n}\) conjecturally admit a
motivic interpretation.  

The Tate conjecture (Conjecture \ref{tatec}) predicts that
\[
 \rho_n = {\rm rank}~\Hom(\Z(-n), h^{2n}(X)), 
\]in the category of motives and \(h^{2n}(X)\) is part of the motive
of \(X\). A related variant is that \(\rho_n\) is the rank of the
Chow group \(CH^n(X)\) of algebraic cycles of codimension \(n\)  on
\(X\). 

The Lichtenbaum-Milne conjecture (Conjecture \ref{lmc}) expresses
\(a_X(n)\)   as an Euler-characteristic of \'etale motivic
cohomology  \(H^*(X,\Z_X(n))\); i.e.,  
the cohomology of the (\'etale) motivic complexes
\(\Z_X(n)\) of S.~Lichtenbaum;  \(\Z_X(0)\) is the constant sheaf \(\Z\)
and \(\Z_X(1)=\G_m[-1]\) is the sheaf \(\G_m\) in degree one.
Their conjecture is known for \(n=0\) (unconditionally) and for \(n=1\)
(modulo the Tate conjecture for divisors on \(X\));
cf. \cite{sl1984, mi5}. For \(n=0,1\), it takes the form 
\[
a_X(0) = \pm \chi(X, \Z), \quad \quad a_X(1) = \pm\frac{\chi(X, \ga)}{\chi(X,\gm)} = \pm\frac{q^{\chi(X,
    \mathcal O_X)}}{\chi(X,\gm)}. 
\]
Lichtenbaum \cite{sl2003} has provided another elegant interpretation
of \(a_X(0)\) using his Weil-\'etale  topology (cf. Theorem \ref{sl}). 

Let us now turn to the special value at \(s=\half\) or
\(t= 1/{\sqrt{q}}\). One can ask for a motivic
description of  

\(\bu\) the order of
vanishing \(\rho_X:= {\rm ord}_{s=\half}\zeta(X,s)\) and 

\(\bu\) the corresponding special value \(c_X\) at \(s=\half\), viz.,
  \(c_X:= {\rm lim}_{t \to 1/\sqrt{q}}  (1- 
  \sqrt{q}t)^{-\rho_X} Z(X,t)\).  

The main result of this paper provides such a description, under the
condition that \(q=p^{2f}\), i.e., that \(\F_{p^2}
\subset\F_q\).  We note that \(c_X\) may not be rational, if the
condition on \(q\) is dropped.  

Our paper is an exploration,
using the methods of   \cite{sl1984, mi5, sl2003,
  mf1}, of a  beautiful \emph{suggestion} of Yuri Manin that ``a
certain supersingular elliptic curve \(E\) might be useful  in finding
the expression for \(\Z(\half)\)'', because its one-motive \(h_1(E)\) is the
square root of Tate's motive in the same sense as the Dirac operator
is the square root of a Laplacian of a spin manifold.  The rank of \(h_1(E)\) is two
whereas \(\Z(1)\) has rank one. Spin structures on a curve 
involve a square root of the
canonical bundle, the dualizing object or the orientation sheaf. The
\(\ell\)-adic analogue  \(\Z_{\ell}(1)\) of the
dualizing sheaf, on a curve, comes from \(\Z(1)\). This suggests
that a motivic interpretation of \(\zeta(\half)\) should involve a
spin structure on the ``curve'' Spec~\(\Z\). 
We shall refer to \(\Z(\half)\) as
the Manin motive.  Going to 
the double cover Spec~\(\F_{p^2}\) of Spec~\(\F_p\) seems essential in
obtaining the square root \(\Z(\half)\) of \(\Z(1)\): the elliptic
curve \(E\) exists only over \(\F_{p^2}\).

Our starting point was the
observation that  \(E\)
provides the required motivic
interpretation of \(\rho_X\) :  
\[2 \rho_{X} = {\rm rank}~\Hom_{\F_q}(E,Pic(X)),
\](Lemma \ref{ord}) where the abelian variety \(Pic(X)\) is the Picard variety of
\(X\).
This suggests that the \'etale sheaf \(E\)  on \(X\) may provide an
 interpretation of \(c_X\).    

Our main theorem  is an analogue of  Lichtenbaum's result (Theorem
\ref{sl}):
\begin{thm}\label{konege} Let \(X\) be a smooth projective variety over
  \(\F_q\) with \(q=p^{2f}\) and \(E\) as before. Write \(\theta\) for
  the generator of the Weil-\'etale cohomology group \(H^1_W({\rm
    Spec}~\F_q, \Z) ={\Hom}(\Z, \Z) =\Z\); cf. \cite[\S8]{sl2003}.

{\rm (i)} The Weil-\'etale cohomology groups \(H^i_W(X,E)\) are
 finitely generated.

{\rm (ii)} \(H^i_W(X,E) =0\) for \(i > 1+ 2~{\rm dim}~X\). 

{\rm  (iii)} The alternating sum \(\Sigma(-1)^i r_i\) of the ranks \(r_i\) of
\(H^i_W(X,E)\) is zero.  

{\rm
(iv)} The secondary Euler-characteristic \(\Sigma(-1)^i ir_i\) is
equal to \(-2~{\rm  ord}_{s=\half}\zeta(X,s) = -2\rho_X\). In fact,  

\(r_0 =r_1 =
2\rho_X\) and \(r_i =0\) for \(i>1\).  

{\rm
(v)} The cohomology groups \(h^i\) of the complex
\[
H^0_W(X,E) \xrightarrow{\cup \theta} H^1_W(X,E) \xrightarrow{\cup
  \theta} H^2_W(X,E) \to \cdots
\]are finite;  the special value at \(s=\half\), viz., \(c_X:= {\rm lim}_{t \to 1/\sqrt{q}}  (1- 
  \sqrt{q}t)^{-\rho_X} Z(X,t)\) is given by  
\[
c_X^2 = \frac{q^{\chi(X,\mathcal O_X)}}{\chi(X,E)}
\]where 
\[
\chi( X,E) = \prod[h^i]^{(-1)^{i}}.
\]
\end{thm} 
Tate's theorem
\cite{Tatendo} on abelian varieties is crucial here; since this result
of Tate is a special case (of divisors on abelian varieties) of his
conjecture,  the case of \(s=\half\) is
intermediate between the case of \(s=0\) and \(s=1\). As \(h_1(E)\) has
rank two, it is \(2\rho_X\) and
\(c_X^2\) that arise rather than \(\rho_X\) and \(c_X\). 
Note that
the Weil-\'etale motivic cohomology groups \(H^*_{W}(X, \Z_X(n))\) are
known to be finitely generated only for \(n=0\) \cite[\S8]{sl2003};  the case \(n\neq0\)
requires the Tate
conjecture \cite{gei}.  

The proof of the main theorem in \S \ref{emc} depends on the results in
\S\ref{patic}. As usual, the main difficulties lie in the \(p\)-part. This
involves the nontrivial computation of the cohomology of the 
flat group scheme \(_{p^n}E\) on \(X\); we do this by using the de
Rham-Witt complex of \(X\); our approach was suggested by
\cite[\S13]{mif}.  The 
special values at \(s=\half\)  for the L-functions of curves 
and motives are treated first in \S  \ref{curva}, \ref {intmotive},
using the work of Milne \cite{mi1, mi2, mf1},   after some preliminaries in
\S\ref{prelim}. The reader may be amused by the results in  \S\ref{bsd}:  the square 
root of the order of a Tate-Shafarevich group turns up in the
description of the special value at \(s=\half\). 

 Even though \(\zeta(X,1/3)\)
will be a nonzero rational  number when \(q=p^{3f}\), one  does not
have \(\Z(\frac{1}{3})\)  over such a field; similar comments apply to \(1/4,
1/5, \cdots\).  The
case of \(\Z(\half)\)  is special and its existence ultimately has its origins in the
structure of Weil \(q\)-numbers \cite{mi6}. We refer to C.~Deninger
\cite[\S7]{Den} and Manin \cite{ma} for ideas on the
exotic Tate motives \(\C(s)\) for \(s\in\C-\Z\).  

Supersingular elliptic curves are important examples of
motives \cite[Theorem~2.41]{mi6}. For instance, as indicated by J.-P. Serre \cite{ag}, their
existence implies the non-neutrality of  the Tannakian category 
\(\mathcal M\)  of motives over a finite field (i.e. that \(\mathcal
M\) does not have a fibre 
functor to \(\Q\)-vector spaces).

Finally, elliptic curves or abelian varieties cannot provide
\(\Z(\half)\) over Spec~\(\Z\): the Hodge numbers are incompatible.
Since elliptic curves with noncommutative endomorphism 
 rings provide \(\Z(\half)\) over a finite field, it may be that the
 arithmetic theory  of non-commutative tori (initiated by Manin) holds
 the key to the definition of \(\Z(\half)\) over Spec~\(\Z\).\medskip

\noindent {\bf Notations.} 
Let  \(\F_q\)  be a finite field of characteristic \(p >0\) and let 
\(\F\) be an algebraic closure of \(\F_q\).  For any scheme \(V\) over
\(\F_q\), we write \(\bar{V}\) for its base change to \({\F}\). Write
\(\Gamma\) for the Galois group \({\rm Gal}(\F/{\F_q})\) and
\(\Gamma_0 \cong \Z\) for the subgroup (the Weil group) generated by the Frobenius
automorphism \(\gamma\): \(x \mapsto x^q\). The additive and
multiplicative valuations of \(\Q\) are normalized so that \(|p|_p =
1/p\) and \({\rm ord}_p(p)=1\). 
We write \(H\) (resp. \(H_{fl}\), \(H_{W}\)) for the \'etale
(resp. flat, Weil-\'etale \cite{sl2003}) cohomology groups. Finally, \([S]\) denotes
the order of a finite set \(S\). 

From  \S \ref{curva} onwards, we assume that \(q\) is a square, i.e.,
that \(\F_{p^2} \subset\F_q\).

\section{Preliminaries}\label{prelim} Here we recall some relevant facts and
conjectures. We begin with an example.

\subsection{Example}  Let us consider the scheme \(X = {\rm
  Spec}~\F_q\). The value \((1-q^n)^{-1}\) of \(\zeta(X, s)
=(1-q^{-s})^{-1}\) at any negative integer \(-n\) has many  interpretations:
\[
-\frac{1}{\zeta( X,-n)} = (q^n-1) = [K_{2n-1}(\F_q)] =
[\Ext^1_{\mathcal A}(\Z, \Z(n))] = [\G_m(\F_{q^n})] = [T_n(\F_q)] 
\]
where \(K_{2n-1}(\F_q)\) is the Quillen K-group of \(\F_q\), \(\mathcal A\) is the
category of effective integral motives over \(\F_q\) of
\cite{mf1}, \(T_n\) is the torus over \(\F_q\) obtained from \(\G_m\)
over \(\F_{q^n}\) via Weil restriction of scalars.

The value of
\(\zeta(X,s)\) at \(s=-\half\) is \((1-q^{\half})^{-1}\) which is not
rational unless \(q\) is an even power of \(p\). So assume that \(q=
p^{2f}\). Even then, 
one cannot interpret \(p^f-1\) as the order of \(\F_q\)-points of any
group (for varying \(f\)). But all is not lost. Consider 
\(\zeta(X,-\half)^{-2} = (p^f-1)^2\) which could be the order of
\(\F_q\)-rational points of an elliptic curve \(E\). In fact, one has  
\[\frac{1}{\zeta(X,- \half)^{2}} =   (p^f-1)^2 = [E(\F_q)] =[\Ext^1_{\mathcal
  A}(\Z, h_1(E))]  
\]for an elliptic curve \(E\) with Frobenius eigenvalues
 \((p^f, p^f)\).  Such an \(E\) is supersingular. 

As for \(s=\half\), we find that 
\[
\zeta(X,\half)^2 = \frac{q}{(1-p^f)^2} = \frac{q^{\chi(\mathcal O_X)}}{[E(\F_q)]}
\]where \(\chi(\mathcal O_X)\) is the
Euler-characteristic of \(\mathcal O_X\).  
Note that the numerator \(q^m\) of \(\zeta(X,m)=\frac{q^m}{q^m-1}\)
for any positive integer \(m\) is an Euler-characteristic by 
Conjecture \ref{lmc}; cf. \cite[Thm.~7.2]{mi5}.
\qed

\subsection{Supersingular elliptic curves and
  \(\Z(\half)\)}\label{zalf} We indicate the relations of supersingular elliptic
curves to the motive \(\Z(\half)\) by quoting  the email message
(dated 10 August 2004) of Milne:

``A.~Grothendieck (and P.~Deligne) knew already in the 1960's that a
supersingular elliptic curve \(E\) over \({\F}\) provides a square
root of the Tate motive in the following precise sense: after
extending the coefficient field to a field \(k\supset \Q\) splitting
\(\End(E)\), \(h_1(E)\) decomposes into \(\Q(\half)\oplus\Q(\half)\)
where \(\Q(\half)^{\otimes2} =\Q(1)\). This is exploited throughout
\cite{mi6}. For example, as for any category of motives, there is an
exact sequence 
\[
0 \to \gm \xrightarrow{\omega } P \xrightarrow{t} \gm \to 0
\]such that \(t\circ\omega=2\) where \(P\) is the kernel of the Tannakian
groupoid. The map \(\omega\) is given by the weight gradation and \(t\) by
the Tate motive. From \(E\) one gets a homomorphism \(P
\xrightarrow{e} \gm\) such that \(e\circ\omega=1\), and so \(P\) decomposes into
\(P = P^0 \times w(\gm)\); cf. \cite[2.41]{mi6}. 

Over a finite field \(\F_q\), supersingular elliptic curves come in
three types:

(a) eigenvalues \(\pm\sqrt{-q}\);

(b) \(q=p^{2f}\); eigenvalues \(-p^f, -p^f\);

(c)  \(q=p^{2f}\); eigenvalues \(p^f, p^f\).

This follows from Tate's theorem \cite[Thm.~1(a)]{Tatendo} but was
probably known to Deuring. The ones in (c) play the same role over
\(\F_{p^{2f}}\) as \(E\) over \(\F\).''\medskip 
 
\nd {\em Remark.} The elliptic curves in (c) are all \(\F_q\)-isogenous. Let \(E\)   be
an elliptic curve over \(\F_{p^{2f}}\) of type (c). It is actually defined over
\(\F_{p^2}\) and its endomorphism ring \(\End_{\F_q}(E) = \End_{\F}(E)\)
 is a maximal order in a quaternion algebra over
\(\Q\), ramified only at \(p\) and \(\infty\); cf. p. 528 and Theorems
4.1, 4.2 of \cite{Wat}.\medskip  

\begin{quote}

{\em We fix an elliptic
curve of type (c) over \(\F_{p^{2}}\) and we refer to it
throughout as \(E\).}\medskip  

\end{quote}

\nd  {\em Remark.}
An abelian variety \(A\) over \(\F_{p^{2f}}\)
whose Frobenius spectrum is pure with support  \(p^f\) is, by Tate's
theorem \cite[Thm.~1(a)]{Tatendo} (cf. \cite[pp.526-27]{Wat}),
isogenous to a power of \(E\).  

\subsection{Conjectures of Tate and Lichtenbaum-Milne}  
Let \(X\) be a smooth projective (geometrically connected) variety over
\(\F_q\). Let \(\zeta(X,s)\) be its zeta function. The associated
function \(Z(X,t)\), defined by \(\zeta(X,s) = Z(X, q^{-s})\), is known
to be a rational function of \(t\) with integer coefficients. The
special value of \(\zeta(X,s)\) at \(s=n\) is the special
value \(a_X(n)\) of \(Z(X,t)\) at \(t=q^{-n}\), up to powers of log\((q)\). 

A motivic  description of \(a_X(n)\) requires a knowledge of
the poles of \(\zeta(X,s)\) given by 

\begin{conj}\label{tatec}{\rm (Tate)} \(T(X,r):\) For any integer \(r \ge 0\), the dimension of the
  subspace of \(H^{2r}(\bar{X}, \Q_{\ell}(r))\) generated by algebraic
  cycles is equal to the order \(\rho_r\) of the pole of \(\zeta(X,s)\) at
  \(s=r\). 
\end{conj}
This conjecture provides a motivic description of \(\rho_r\), i.e.,  \(\rho_r = {\rm
  rank}~\Hom(\Z(-r), h^{2r}(X))\) in the category of motives over
\(\F_q\). While Conjecture \ref{tatec} is unknown in general, even for
divisors \(r=1\), it does hold for divisors on abelian varieties
\cite[Thm.~1(a)]{Tatendo}.

In the case of a curve \(X\), it is well known that
\(a_X(0), a_X(1)\) are related to the class number \(h\) of \(X\)
\cite[p.~191]{sl1984}.

The description of \(a_X(1)\)  for
a surface \(X\) can be considered,  following M.~Artin and J.~Tate
\cite{Ta}, as a geometric analogue of the 
Birch-Swinnerton-Dyer conjecture. 
Inspired by \cite{Ta}, Lichtenbaum \cite{sl1984} and Milne
\cite{mi5} have conjectured 
(Conjecture \ref{lmc})  a complete description (including the
\(p\)-part) of \(a_X(n)\) for any \(X\) as an Euler-characteristic in \'etale motivic
cohomology  \(H^*(X,\Z_X(n))\).  Namely, Lichtenbaum (resp. Milne) has
defined 
a non-zero rational number \(\chi(X,\Z(r))\) (resp. an integer
\(\chi(X, \mathcal O, r)\)) via motivic cohomology \(H^*(X, \Z_X(r))\)
(resp. \(\Sigma_{i=0}^{i=r}(r-i)\chi(X, \Omega^i)\)) of \(X\); we refer to \cite[\S0]{mi5}
for the precise definitions. These are related to \(a_X(n)\)  via the 

\begin{conj}\label{lmc}{\rm (Lichtenbaum-Milne)} \cite[0.1]{mi5} The special
  value \(a_X(r)\) of \(Z(X,t)\) at \(t=q^{-r}\) can be given as 
an Euler-characteristic of motivic cohomology. Namely, one has   
\[
\tag{LM(X,r)} Z(X,t) \sim \pm \chi(X,\Z(r)). q^{\chi(X,\mathcal O,
  r)}  (1-q^{r}t)^{-\rho_r}, \quad {\rm as}~ t \to q^{-r} 
\]where the terms on the right are defined  in \cite[Conj. 0.1]{mi5} 
\end{conj}

This conjecture generalizes \cite{sl1984} that of Artin-Tate \cite{Ta} for the case
of \(r=1\)  and \(X\) a surface. We recall one result about Conjecture
\ref{lmc} and refer to \cite{sl1984, mi4, mi5, mf1}  for other 
results.
\begin{thm}  \cite{sl1984, mi4} {\rm (a)}   Conjecture \ref{lmc}  
is true for \(r=0\).  

{\rm  (b)} If \(T(X,1)\) holds, then  the
terms in {\rm (LM(X,1))} are finite and {\rm (LM(X,1))} is true.   
\end{thm} 

Part (b) generalises the result of Artin-Tate \cite{Ta} for surfaces, and includes the \(p\)-part  \cite{mi3}.   

\subsection{Weil-\'etale cohomology}\label{weit} As before, \(X\) is a smooth
projective geometrically connected variety.  Lichtenbaum
\cite{sl2003} has given a beautiful interpretation, via his
Weil-\'etale  cohomology groups \(H^*_W(X,\Z)\), of the behaviour 
of \(Z(X,t)\)  at \(t =1\).  Namely, cup-product with the 
generator \(\theta\) of \(H^1_W({\rm Spec}~\F_q, \Z) ={\Hom}(\Z, \Z) =\Z\) gives a
complex \((H^{\bullet}_W(X,\Z), \theta)\) of finitely generated
abelian groups whose cohomology groups
\(h^*_W(X)\) are finite; write \(\chi_W(X)=
\frac{[h^0_W(X)]\cdots}{[ h^1_W(X)]\cdots}\) for the alternating
product of the orders \([h^*_W(X)]\) of \(h^*_W(X)\). Then, Lichtenbaum's
result\footnote{His results include the case of  arbitrary curves and smooth surfaces.}  on
the behaviour of \(Z(X,t)\) at \(t=1\) is as follows: 

\begin{thm}\label{sl} {\rm (Lichtenbaum)} 

{\rm (i)} the usual Euler characteristic \(\Sigma(-1)^i r_i\) of
\(X\) is zero; here \(r_i = {\rm rank}~H^i_W(X,\Z)\).

{\rm (ii)} the  order of
the zero of \(Z(X,t)\) at \(t=1\)  is given by the \emph{secondary
  Euler characteristic} \(\Sigma(-1)^i i . r_i\). 

{\rm (In our case, this
order is minus one.) } 

{\rm (iii)} the special value \(a_X(0):={\rm lim}_{t \to 1}Z(X,t)
(1-t)\)  is the \emph{multiplicative} Euler characteristic: 
\[a_X(0) =\pm\chi_W(X).
\]\end{thm}

\subsection{Behaviour at \(s=\half\): basic invariants} Our aim in
this paper is a motivic description, when \(q=p^{2f}\), of the following arithmetic
invariants of \(X\):  
\begin{itemize} 
\item the integer \(\rho_X:=\)  the order of
the zero of \(\zeta(X,s)\) at \(s=\half\). 

\item   the special  value at \(t=q^{-\half}\) of \(Z(X,t)\), viz.,
  \(c_X:= {\rm lim}_{t \to 1/\sqrt{q}}  (1- 
  \sqrt{q}t)^{-\rho_X} Z(X,t)\).
\end{itemize} 

As before, \(c_X\) differs from the special value at \(s=\half\) of
\(\zeta(X,s)\) by factors of \(log~p\).   As \(q=p^{2f}\), one has 
\(c_X \in\Q^*\).  Unlike \(\rho_r\), the integer \(\rho_X\) has an unconditional
motivic description, viz., 

\begin{lem}\label{ord} Let \(q=p^{2f}\) and let \(Pic(X)\) be the
  Picard variety of \(X\). One has  
\[2 \rho_X=  {\rm rank}~\Hom_{\F_q}(E,
Pic(X)).
\]
\end{lem}  

For instance, \(\rho_E =2\) and  \({\rm rank}~\End(E) = 4\). The RHS
is always divisible by \(4\) because \(\Hom_{\F_q}(E,
Pic(X))\) is a module over \(\End(E)\).  Similarly, the LHS is 
divisible by \(4\): the integer \(\rho_X\)  is even because the roots 
\(\alpha \neq \sqrt{q}\) of the factor \(P_1(t)\) of \(Z(X,t)\) come
in pairs \((\alpha, q/{\alpha})\) and the degree of \(P_1(t)\) is
even.  
 
\begin{proof} This follows from Tate's theorem
\cite[Thm.~1(a)]{Tatendo}: by the Weil conjectures, \(\rho_X\) is the
multiplicity of  the root \(\sqrt{q}\) in the factor \(P_1(t)\) of
\(Z(X,t)\) which is the characteristic polynomial of Frobenius on
\(Pic(X)\). For \(A=E\) and \(B=Pic(X)\), the integer \(r(f_A, f_B)\)
defined by Tate \cite[(4), p.~138]{Tatendo}  is \(2\rho_X\).\end{proof}

Thus, the vanishing of \(\zeta(X,s)\) at \(s=\half\) is controlled by
the (non)-ordinarity of \(Pic(X)\). The study of \(c_X\) occupies the
rest of the article. 

\section{Curves}\label{curva} 
In this section, the work of Milne \cite{mi1,mi2} is used to
provide motivic interpretations of \(c_X\) in the case of curves.  Let  \(X\) be a
curve of genus 
\(g\) over \(\F_q\) with \(q=p^{2f}\), i.e., \(\F_{p^2} \subset
\F_q\).  Write \(J\) for its Jacobian. 

\subsection{Extensions and \(c_X\)} 

\begin{thm}\label{curve} (i) \(\rho_X = \half {\rm rank}~\Hom_{\F_q}(E,J)\). 

(ii) the special value at \(s=\half\) of \(\zeta(X,s)\) is given by  
\[c_X^2= q^{\chi(\mathcal O_X)}. \frac{[\Ext^1_{\F_q}(J,E)]}{[E(\F_q)]^2}.D
\]where \(\chi(\mathcal O_X) = 1-g\) is the Euler-characteristic of
\(\mathcal O_X\), \(\Ext^1\)  is computed in the category of group
schemes over \({\rm Spec}~\F_q\), and \(D\) is the discriminant of the pairing 
\[\Hom_{\F_q}(E,J) \times \Hom_{\F_q}(J,E)  \xrightarrow{\circ} \End_{\F_q}(E)
\xrightarrow{trace} \Z. 
\]

(iii) one has \[c_X^2= q^{\chi(\mathcal O_X)}. \frac{[H^1(X,E)]}{[H^0(X,E)_{tors}]^2}.D
\]where \(H^*(X,E)\) is the \'etale cohomology of the sheaf 
defined by \(E\) and \(D\) is as in (ii). 
\end{thm} 

\begin{proof} (i) was proved earlier. 

(ii) When \(g=0\), one easily verifies that 
\[
c_{{\Pe}^1}^2 = \frac{q}{[E(\F_q)]^2}.
\] When
\(g\neq 0\), \(Z(X,t) = P(t). Z({\Pe}^1,t)\) where
\(P(t) = \prod (1-\alpha_it)\), the numerator
of \(Z(X,t)\), is the characteristic polynomial of Frobenius on
\(J\). The result now follows from Lemma (\ref{ord}) and  
\cite[Thm.~3]{mi1}.

(iii) As \(\Ext^1_{\F_q}(J,E) \cong H^1(X,E)\) \cite[Cor.~3]{mi2} and \(E(\F_q)=
H^0(X,E)_{tors}\) \cite[p.~120]{mi2}, this is clear.\end{proof} 

\nd {\em Remark.}  If \(X\) is ordinary, i.e., its Hasse-Witt matrix is
  invertible or the abelian variety \(J\) is ordinary,
then the formula simplifies to 
\[
\zeta(X,\half)^2 = c_X^2 =
q^{1-g}\frac{[\Ext^1_{\F_q}(J,E)]}{[E(\F_q)]^2}.  
\]Since \(J\)  is ordinary, the roots of \(P\) satisfy: \(\alpha_i
\neq \sqrt{q}\); thus 
\(\zeta(X,\half)\neq0\). 
For an ordinary elliptic curve \(X\), a simple computation yields  
\[\pm c_X = 1 -
\frac{[X(\F_q)]}{[E(\F_q)]}.
\]

\subsection{BSD and \(c_X\)}\label{bsd} These formulas for \(c_X^2\) are very reminiscent of
the Birch-Swinnerton-Dyer (BSD) conjecture. In fact, Milne \cite{mi2}
has 
shown \(\Ext^1_{\F_q}(J,E) \cong \Sha (E/K)\) of the constant elliptic curve 
\(E\) over the function field \(K=\F_q(X)\).  Since the order of
\(\Sha (E/K)\) is a square,  \([\Ext^1_{\F_q}(J,E)] = m^2\) for a positive 
integer \(m\). This 
gives a motivic interpretation of \(c_X\) for \(X\) ordinary: 
\[
c_X= \pm p^{(1-g)f}
\frac{m}{[E(\F_q)]}.\]Thus, the square root of the order
of \(\Sha\) is related to special values at \(s=\half\). 

\begin{thm} The special value at \(s=1\) of the L-function
  \(L(E/K,s)\) of \(E\) over \(K\) is:  
\[{\rm lim}_{s\to 1} L(E/K,s). (s-1)^{-2\rho_X} = c_X^2 .({\rm log}~q)^{\rho}
\]where \(\rho  ={\rm rank}~E(K) = 2 \rho_X = {\rm rank}~\Hom_{\F_q}(J,E)\).
\end{thm}

\begin{proof} This follows from \cite[Thm.~3]{mi2} and the first
equation on page 102 of  \cite{mi2}.\end{proof}

\section{Integral motives}\label{intmotive}Some of the results in the
previous section generalize  to effective integral motives \cite{mf1}.  

Recall the category \(\mathcal M^+(\F_q
;\Z)\) of effective integral 
motives over \(\F_q\) and the category
\(\mathcal M(\F_q ;\Z)\) of integral motives (see
\cite{mf1}). Assume that \(q=p^{2f}\). 
The elliptic curve \(E\) defines an effective integral motive
\cite[5.16]{mf1} via its 
\(h^1\) which is still denoted \(E\). Its dual \(E^{\vee} = h_1(E)\) is just \(E\otimes\Z(1)\).  The special value at 
\(s=\half\) of the L-function \(L(M,s)\) of any effective integral
motive \(M\) is related to the order of 
Ext-groups in \(\mathcal M^+(\F_q ;\Z)\):

\begin{prop}\label{intmot}Let \(q=p^{2f}\) and let  \(r\) be the rank
  of \(M\) and \(\rho_M:={\rm ord}_{s=\half}L(M,s)\). One has

{\rm (i)} \(2\rho_M= {\rm rank}~\Hom(E,M)\). 

{\rm (ii)}  
\[q^{r}L(M,s+\half)^2 \sim \frac{[\Ext^1(E, {M})] .
  D({M})}{[\Hom(E,{M})_{tors}].[\Hom({M},E)_{tors}]}.  
(1-q^{-s} ) ^{\rho}    \quad {\rm as} \quad s \to 0,
\]where  \(D({M})\)  is the discriminant of the pairing 
\[\Hom(E,M) \times \Hom({M},E) \xrightarrow{\circ}
\End({E}) \xrightarrow{trace} \Z.
\] 

{\rm (iii)} One has \(L(M,s +\half)^2  = L({M}\otimes E^{\vee},s).\)
\end{prop} 

\begin{proof} If \(L(M,s) = \prod  (1-a_iq^{-s})\), then by \cite[9.2(d)]{mf1},  \[ 
L(M\otimes E^{\vee}, s) = \underset{i,j}{\prod}(1-
\frac{a_i}{b_j}q^{-s}) = \underset{i}{\prod}(1- a_i q^{-s-\half})^2 
\]because \(L(E,s) =  \prod(1-b_jq^{-s}) =
(1-q^{\half}q^{-s})^2\). This proves (iii).

Now it is clear that (i) and (ii) follow from \cite[Thm
10.1]{mf1}. \end{proof}  
 
It is possible to prove  a  Weil-\'etale variant \cite{sl2003} of the
above result using \cite[Thm. 10.5]{mf1}. 

\nd {\em Remark.}  (i)  One has 
\[
L(\Z, \half)^2 = (1-q^{-\half})^2 = \frac{[E(\F_q)]}{q}
\]By the proposition, \(qL(\Z, \half)^2\) is  \([\Ext^1(E, \Z)]\)
which, by
\cite[8.7]{mf1}, is \([E(\F_q)]\).

(ii) Let \(A\) be an ordinary abelian variety of dimension \(d\). One
computes \(L(h_1(A),\half) \neq 0\) to be \cite{mi1} 
\[L(h_1(A), \half)^2 = q^{-d} [\Ext^1(A,E)] 
\]where the {\rm Ext}-group is computed in the category of group
schemes over \(\F_q\). 

(iii)   For any motive \(M\) over a global (or finite) field and for 
any \(n \in \Z\), we write \(M(n)\) for its (Tate) twist by
\(\Z(n)\). On the level of L-functions, 
twisting corresponds to translations on the \(s\)-axis: \(L(M,s+n) =
L(M(n),s)\).  The above proposition computes special values
of \(L(M,s)\) at other half-integers as well.\qed

\section{The \(p\)-adic Tate module of \(E\)}\label{patic}
As before, \(X\)  is a smooth projective geometrically
connected variety over Spec~\(\F_q\) with \(q=p^{2f}\), and \(E\) is
our fixed elliptic curve. 

The main result (Theorem \ref{bures}) of  this section relates \(c_X\)
to the cohomology groups \(H^*_{fl}({X},~_{p^n}E)\) of the finite
flat group scheme \(_{p^n}E\) on \(X\).  This result is the technical core
of the paper.

We shall freely use the
results and methods of \cite{mi4}. We recall from \cite[\S1]{mi4} the
category \(Pf\) of perfect affine schemes over \(\F_q\), endowed with
the \'etale topology.  The computation of \(H^*_{fl}(\bar{X},~_{p^n}E)\) here was inspired by  \cite[\S13 (a)]{mif}.

For any perfect field \(k\) of positive characteristic, we write
\(W(k)\) for the Witt ring of \(k\).  Recall the Dieudonne ring \(A =
  W(\F_q)_{\sigma}[F,V]\); put \(\bar{A} =A\otimes_{W(\F_q)}W(\F)\).
  Note that the relation \(FV = p  =VF\) holds in \(A\) and \(\bar{A}\).

\begin{lem}\label{tamu} (i) For \(n \ge 1\) and \(i \ge 0\), one has an exact sequence  
\[ \cdots H^i_{fl}(\bar{X},~ _{p^n}E) \to H^i(\bar{X}, W_{2n}\mathcal O)
\xrightarrow{F-V} H^i(\bar{X}, W_{2n}\mathcal O) \to \cdots  ,
\]where \(W_{2n}\mathcal O\) is the Witt vector sheaf of order \(2n\) on
\(\bar{X}\).

(ii) The presheaf \(T \mapsto
H^i({X}_T,~ _{p^n}E)\) on \(Pf\) is represented by an affine perfect
group scheme 
\cite[pp.303-4]{mi4}. \end{lem} 
\begin{proof} (i) This follows from the exact sequence of flat sheaves
\[
0 \to~ _{p^n}E \to W_{2n}\mathcal O \xrightarrow{F-V} W_{2n}\mathcal O
\to 0\] on \(\bar{X}\), a consequence of the fact that the
Dieudonne module of 
  the p-divisible group \(E(p^{\infty})\)  (resp. \(W_{2n}\)) is
  \(E':= \bar{A}/(F-V)\) (resp. \(\bar{A}/(V^{2n})\)). The cohomology
  of \(W_{2n}\mathcal O\) is the same in the \'etale and flat
  topologies. 

(ii) Follows from (i) in a standard manner, cf.  \cite[Lemma
1.8]{mi4}, because the presheaf \(T\mapsto H^i(X_T, W_{2n}\mathcal
O)\) is representable.\end{proof}

As in \cite [p. 322]{mi4}, the long exact sequence of Lemma \ref{tamu}
(i) can be regarded as one of affine 
perfect group schemes. Write \(\mathcal H^i(X, W\mathcal
O):=\varprojlim \mathcal H^i(X, W_{2n}\mathcal O)\) and \(\mathcal
H^i(X, ~_{p^n}E)\) for the perfect group
scheme in (ii).  Write \(b_i\) for the
dimension of the perfect pro-algebraic group scheme \(\mathcal H^i(X,
T_pE):= \varprojlim \mathcal H^i(X, ~_{p^n}E)\). Thus, \(b_i =\infty\)
unless the neutral component \(\mathcal H^i(X, T_pE)^0\) of \(\mathcal
H^i(X,T_pE)\) is algebraic, in which case \(b_i\) is the number of
copies of \(\ga^{pf}\) occuring as quotients in a composition series for \(\mathcal H^i(X,
T_pE)^0\).  

Write \(W\) for \(W(\F)\). There is an exact sequence \cite[ p.~321]{mi4}
\begin{equation}\label{dodo}
0 \to H^i(\bar{X}, W\mathcal O)_t \to H^i(\bar{X}, W\mathcal O) \to B_i \to 0,
\end{equation}of \(W[[V]]\)-modules in which \(B_i\) is a free \(W\)-module of
finite rank and \( H^i(\bar{X}, W\mathcal O)_t\) is a finitely
generated \((W/{p^nW})[[V]]\)-module for some \(n\). Write \(d_i\) for
the length of \(H^i(\bar{X}, W\mathcal O)_t \otimes_{W[[V]]}W((V))\)
as a \(W((V))\)-module. The integer \(d_i\) is equal to the number of
copies of \(\F[[V]]\) occuring in a composition series for
\(H^i(\bar{X}, W\mathcal O)_t\). 

\begin{prop}\label{nice}With notations as above, \(b_i\) is finite for all \(i\) and
  \(\Sigma(-1)^ib_i = -\Sigma(-1)^id_i\).\end{prop} 

\begin{proof} We follow the proof of \cite[Prop.~3.1]{mi4} to which we
  refer for the properties of the integer
  \(\chi(\alpha):= {\rm dim}\Ker(\alpha) - {\rm dim}\Coker(\alpha)\)
  attached to a morphism \(\alpha\) of group schemes. 

 We first
  compute \(\chi(F-V)\) on \(H^i(\bar{X}, W\mathcal O)\), in the
  terminology of \cite[Prop.~3.1]{mi4}. By \cite[3.2 (b)]{mi4}, it is
  possible to neglect finitely generated torsion
  \(W\)-modules which arise as subquotients of \(H^i(\bar{X},
  W\mathcal O)\). As \(F=0\) on \(\F[[V]]\), one
  has \(\chi(F-V|\F[[V]]) = \chi(-V|\F[[V]]) = -1\). So
  \(\chi(F-V|H^i(\bar{X}, W\mathcal O)_t) = -d_i\). 

We claim that  \(\chi(F-V|H^i(\bar{X}, W\mathcal O)) =
\chi(F-V|H^i(\bar{X}, W\mathcal O)_t)\).  Namely, we claim that 
\(\chi(F-V|B_i) =0\). This is a consequence of the semisimplicity of
the category of
isocrystals over \(\F\): namely, the kernel (resp. cokernel) of
\(F-V\) on \(B_i\) is \(\Hom_{\F}(E',B_i)\)  (resp. \(\Ext^1_{\F}(E',B_i)\)) in the
category of Dieudonne modules over \(\F\). Here \(E':= \bar{A}/(F-V)\)
denotes  the Dieudonne module of the \(p\)-divisible group
\(E(p)\). As \(\Hom_{\F}(E',B_i)\) is a 
finitely generated free \(\Z_p\)-module and \(\Ext^1_{\F}(E',B_i)\) a \(p\)-primary
torsion group whose \(p\)-torsion is finite, one has \(\chi(F-V|B_i)=0\). 

Writing \(K_i\) and \(C_i\) for the kernel and cokernel of \(F-V\) on
\( \mathcal H^i(\bar{X}, W\mathcal O)\), one has an exact sequence 
\begin{equation}\label{bala}
0 \to C_{i-1} \to \mathcal H^i(\bar{X}, T_pE) \to K_i \to 0. 
\end{equation}
This, as in \cite[p. 323]{mi4}, proves the proposition.\end{proof} 

Any group \(G\)  in \eqref{bala} is an extension of an
\'etale group \(G^{et}\) by a connected group \(G^0\). Modulo finite
groups, one has \(K_i^{et}(\F) =\Hom_{\F} (E',B_i)\) and \(C_i^{et}(\F) =
\Ext^1_{\F}  (E',B_i)\), and an exact sequence
\[
0 \to\Ext^1_{\F} (E',B_{i-1}) \to \mathcal H^i(\bar{X}, T_pE)^{et}(\F) \to
\Hom_{\F} (E', B_i) \to 0.
\]Note \(H^i(\bar{X}, T_pE) = \mathcal H^i(\bar{X}, T_pE)(\F)\). Set
\(H^*(X, T_pE):=\varprojlim H^i_{fl}(X,~_{p^n}E)\). 

\begin{lem} There is an exact sequence
\[
0 \to H^{i-1}(\bar{X}, T_pE)_{\Gamma} \to H^i(X, T_pE) \to
H^i(\bar{X}, T_pE)^{\Gamma} \to 0.
\]
\end{lem} 

\begin{proof} Easy adaptation of the proof of 
\cite[Lemma~3.4]{mi5}.\end{proof} 
 
\subsection{The action of \(\gamma\) on \(H^*(\bar{X}, T_pE)\)} We now
come to the main step in relating the \(p\)-adic Tate module of \(E\)
to the special value \(c_X\). 

Consider the map \(\alpha_i:H^i({X}, T_pE) \to H^{i+1}({X},
T_pE) \) defined by the following commutative diagram  
\[
\begin{CD} 
 H^i({X}, T_pE) @>{\alpha_i}>> H^{i+1}({X}, T_pE)\\
@VVV @AAA\\
H^i(\bar{X}, T_pE)^{\Gamma} @>>> H^i(\bar{X}, T_pE)_{\Gamma};
\end{CD}
\]the vertical maps arise from the previous lemma. As the lower map is
cup-product with the canonical generator \(\theta_p \in H^1(\Gamma,
W)\) and \(\theta_p^2 =0\), 
\[
M^{\bu}: \quad \cdots \to H^i({X}, T_pE) \xrightarrow{\alpha_i} H^{i+1}({X},
T_pE) \to \cdots
\]is a complex; cf. \cite[pp. 544-545]{mf1} for details.  Define 
\[
z = \prod [H^i(M^{\bu})]^{(-1)^i}
\]when these numbers are finite. 
Write, as usual, \(Z(X,t) = \prod P_i(t)^{(-1)^{i+1}}\) and \(P_i(t) =
\prod_j(1-\omega_{ij}t)\) the characteristic polynomial of \(\gamma\)
on  \(\ell\)-adic cohomology \(H^i(\bar{X}, \Z_{\ell})\). 
Our next result is a variant of \cite[Theorem~9.6]{mf1} (cf. \S\ref{weit}).
\begin{thm}\label{bures} (i) \(H^*(X, T_pE)\) are finitely
  generated \(\Z_p\)-modules.  

(ii) the usual Euler characteristic \(\Sigma(-1)^i~{\rm rank}~H^i(X,
T_pE)\) is zero. 

(iii) the secondary Euler characteristic \(\Sigma(-1)^i i . {\rm
  rank}~H^i(X, T_pE)\) is \(2\rho_X\). 

(iv) the cohomology groups \(H^i(M^{\bu})\) are
  finite. 

(v) the value of \(z\)  is given by 
\[  z = \prod_i([H^i(M^{\bu})])^{(-1)^{i}}=   \Bigl\lvert \frac{q^{\chi(X, \mathcal O_X)}}{c_X^2}\Bigr\rvert_p.
\]

(vi) \(H^i(X, T_pE)\) is finite for \(i\neq1,2\). 
\end{thm}

\begin{proof} One can write \(z\)  as a product \(z=z_0z_{et}\)
  corresponding to the decomposition   
\[
0 \to  \mathcal H^i(\bar{X}, T_pE)^0 \to  \mathcal H^i(\bar{X}, T_pE)
\to  \mathcal H^i(\bar{X}, T_pE)^{et} \to 0.
\]of \( \mathcal H^i(\bar{X}, T_pE)\) into its connected and \'etale
parts.  
For each of the groups in the exact sequence 
\[
0 \to C_{i-1}^0 \to   \mathcal H^i(\bar{X}, T_pE)^0 \to K_i^0 \to 0
\]of connected group schemes, the map \(\gamma-1\) 
on the \(\F\)-points is surjective because it is an \'etale
endomorphism of a connected group. Therefore, we obtain that  \(
\mathcal H^i(\bar{X}, T_pE)^0_{\Gamma} =0\), \([ (\mathcal
H^i(\bar{X}, T_pE)^0)^{\Gamma}] = q^{b_i}\), and 
\[z_0 = 
q^{\Sigma(-1)^ib_i} = \Bigl\lvert q^{-\Sigma(-1)^ib_i}\Bigr\rvert_p . 
\]

The finite \(W\)-torsion in \(H^*(\bar{X}, W\mathcal
O)\) does not contribute to \(z\) \cite[pp.~80-81]{mi1}. Thus, we may
assume  that \(K_i^{et}(\F) =\Hom_{\F} (E',B_i)\) and \(C_i^{et}(\F) =
\Ext^1_{\F}  (E',B_i)\), and that there is an exact sequence
\[
0 \to\Ext^1_{\F_q} (E',B_{i-1}) \to H^i({X}, T_pE)^{et} \to
\Hom_{\F_q} (E', B_i) \to 0. \] Now the first term is finite and the last
term is a finitely generated \(\Z_p\)-module.  As \(H^i({X},
T_pE)^{et}\) and \(H^i({X}, T_pE)\) differ only by a finite group of
order \(q^{b_i}\), this proves (i). 

The \'etale part of \(M^{\bu}\) sits in the
commutative diagram
\[
\begin{CD}
{} @.  \Hom_{\F_q}(E', B_1)  @>{\beta_1}>> \Ext^1_{\F_q}(E',B_1) @.
\Hom_{\F_q}(E', B_3)\cdots\\   
@. @AAA @VVV @AAA\\
 H^0({X}, T_pE)^{et} @>{\alpha_0}>> H^{1}({X}, T_pE)^{et}
 @>{\alpha_1}>> H^{2}({X}, T_pE)^{et} @>{\alpha_2}>>
 H^{3}({X}, T_pE)^{et}\cdots\\ 
@VVV @AAA @VVV @AAA\\
\Hom_{\F_q}(E', B_0)  @>{\beta_0}>> \Ext^1_{\F_q}(E',B_0)
@. \Hom_{\F_q}(E', B_2)  
@>{\beta_2}>> \Ext^1_{\F_q}(E',B_2) \cdots\\
\end{CD}
\]where the maps \(\beta_i\) are the ones in
\cite[Lemma~4]{mi1} for the \(p\)-divisible groups \(G_i\) (whose 
dimension we denote by \(g_i\))  whose {\em dual} corresponds to \(B_i\), and
\(H=E(p)\). We shall prove that \(z_{et}\) is defined and calculate
its value by appealing to \cite[Lemma~4]{mi1}. Let us recall two
relevant facts for this purpose. 

\(\bu\) \cite[Remark 5.5]{mi4} the characteristic
polynomial \(P_i(t)=\prod_j(1-\omega_{ij}t):= {\rm det}(1-\gamma t)\) of \(\gamma\) on
\(H^i(\bar{X}, \Q_{\ell})\) is the same as the characteristic
polynomial of \(\gamma\) on the crystalline
cohomology \(H^i_{crys}(\bar{X})\otimes\Q_p\) of \(\bar{X}\).

\(\bu\) \cite[3.5.3, p. 616]{il}  the characteristic polynomial of
\(\gamma\) on \(H^i(\bar{X}, W\mathcal O)\otimes\Q_p\) (by
(\ref{dodo}), this is the same as on \(B_i\)) is
\(\prod_j(1-\omega _{ij})\) where the product is over all \(\omega_{ij}\) with \({\rm
  ord}_q(\omega_{ij}) <1\). 

We can now apply \cite[Lemma~4]{mi1} to obtain that \(z_{et}\) is
defined (which proves (iv)) and given by:\footnote{In the formula for 
  \(z(g)\) in   \cite[Lemma~4]{mi1}, the  exponent of \(q\) should read
  \(d(G^t).d(H)\) - in his notation, this is \(n_2m_1\), as follows from an inspection of the calculation at
  the bottom of p.~81 of \cite{mi1}.}
\[z_{et}= \Bigl\lvert q^{\Sigma (-1)^i g_i}{{\displaystyle\prod_{i}}}(\prod    _{\omega_{ij} \neq
  \sqrt{q}}(1- \frac{\omega_{ij}}{\sqrt{q}}))^{(-1)^i}\Bigr\rvert_p^2
\]where the product is over all \(\omega_{ij}\) with \({\rm
  ord}_q(\omega_{ij}) <1\). Using
\[
\Bigl\lvert 1-\frac{\omega_{ij}}{\sqrt{q}}\Bigr\rvert_p =1 \quad {\rm if}\quad {\rm ord}_q(\omega_{ij}) >
\half,
\]this condition \({\rm
  ord}_q(\omega_{ij}) <1\) may be disregarded in \(z_{et}\). This
proves (iv). 

We now complete the proof of (v); given the formulas for \(z_{et}\) and
\(z_0\), it remains to show that \(\sum (-)^i(g_i
-b_i) = \chi(X,\mathcal O_X)\). As this will use the case \(r=1\) of
\cite[Proposition 4.1]{mi4}, we translate that result into our
context.  First, we observe that our
\(\chi(X,\mathcal O_X)\) is  \(\chi(X, \mathcal
O_X, 1)\) of \cite[Proposition 4.1]{mi4} -- the definition of the
latter is on top of page 325 of \cite{mi4}. Next, our \(d_i\) (defined
just before Proposition \ref{nice}) is Milne's
\(d^i(0)\) (defined on \cite[p. 321]{mi4}, just before Proposition
3.1). 

Now, it remains to interpret the two terms on the right hand side of
\cite[Proposition 4.1]{mi4} with \(r=1\).  The second
term \(\sum_i(-1)^id^i(0)\), which is our
\(\sum_i(-1)^id_i\), can be replaced by \(-\sum(-1)^ib_i\) by
Proposition \ref{nice}. For the first term, we recall that for any
\(p\)-divisible group \(A\) over \(\F_q\), the dimension of the dual group \(A^t\)
is given by the well known (see, for example,  \cite[p. 81]{mi4}; it
also follows from \cite[Thm. 1 (e)]{mi4})
\[
{\rm dim}(A^t) = {\rm height}(A) - {\rm dim}(A).
\]
Now consider the first term in the
right hand side of \cite[Proposition 4.1]{mi4} with \(r=1\). The
numbers \(\lambda_{ij}\) relevant here are those which satisfy
\(\lambda_{ij}\le 1\); by the result \cite[3.5.3, p. 616]{il}
also mentioned earlier, these are exactly the slopes of the roots of the
characteristic polynomial of \(\gamma\) acting on
\(H^*(\bar{X}, W\mathcal O)\otimes\Q_p\) or on the associated \(p\)-divisible
groups \(B_1, B_2, \cdots\) defined earlier. Milne's \(\lambda_{ij}\)
are our \({\rm ord}_q(\omega_{ij})\). The first term can therefore be written as 
\[
\sum_{{\rm ord}_q(\omega_{ij}) \le 1}(-1)^i m_{ij}(1-{\rm ord}_q(\omega_{ij})),
\]where \(m_{ij}\) is the multiplicity of \({\rm ord}_q(\omega_{ij})\). It is now an
easy exercise to see that this sum is \(\sum_i(-1)^i {\rm
  dim}(B_i^t)\). But \(g_i\) is the dimension of the
\(p\)-divisible group whose {\em dual} corresponds to \(B_i\). We
summarise our discussion
\begin{align*}
\chi(X,\mathcal O_X) & = \chi(X, \mathcal O_X, 1)\\
&=\sum_{{\rm ord}_q(\omega_{ij}) \le 1}(-1)^i m_{ij}(1-{\rm ord}_q(\omega_{ij})) + \sum_i(-1)^id^i(0)\\
& =\sum_i(-1)^i ({\rm dim} (B_i^t) + d^i(0))\\
& = \sum_i (-1)^ig_i + \sum_i(-1)^id_i\\
& = \sum_i (-1)^i(g_i - b_i).
\end{align*}This proves (v).

We note that the \(\Z_p\)-ranks of \(\Hom_{\F_q} (E',B_i)\)  and
\(\Ext^1_{\F_q}(E',B_i)\) are equal (the ranks are unchanged by an
isogeny and one reduces to the cases treated in \cite[
pp.81-83]{mi1}). While the first contributes to \(H^i({X},
T_pE)^{et}\), the second contributes to \(H^{i+1}({X}, T_pE)^{et}\),
which proves (ii). 

Now, the rank of  \(\Hom_{\F_q} (E',B_i)\) is
non-zero only if \(\sqrt{q}\) -- a Weil number of weight one -- is a root of the minimal polynomial of Frobenius on
\(B_i\) \cite[p.~81]{mi1}. But a  root \(\omega_{ij}\) of Frobenius on \(B_i\) (or
\(H^i(X, W\mathcal O)\)) is a Weil number of weight \(i\). Thus, the rank
is zero for \(i\neq1\).  This proves (vi). The \(\Z_p\)-rank of \(\Hom_{\F_q}
(E',B_1)\) is, by Tate's theorem \cite{dj}, equal to the \(\Z\)-rank of \(\Hom_{\F_q}(E,
Pic(X))\) which is \(2\rho_X\). As \((-1)^0 0 + (-1) 1\times 2\rho_X + (-1)^2
2 \times2\rho_  X = 2\rho_X\), this proves (iii). \end{proof}

\section{The Weil-\'etale cohomology of an elliptic curve}\label{emc}
As before, \(q=p^{2f}\), \(X\)  is a smooth projective geometrically
connected variety over Spec~\(\F_q\), and \(E\) is our fixed elliptic curve. In
this section, we compute \(H^*(X,E)\) and the Weil-\'etale
\cite{sl2003}) cohomology groups \(H^*_W(X,E)\) and prove Theorem
\ref{konege}.   We write \(A\) for the
Albanese variety of \(X\).

\subsection{Cohomology of the elliptic curve \(E\)}  The basic properties of
\(H^*(X,E)\) are given the  following result whose formulation was
inspired by \cite[Proposition~2.1]{sl1984}.

\begin{thm}\label{manins} (a) The \'etale cohomology of \(E\) is given
  as follows:

(i) \(H^0(X,E)\) is finitely generated and
\[
{\rm rank}~H^0(X,E) = {\rm rank}~\Hom_{\F_q} (A, E) = {\rm rank}~\Hom_{\F_q}(E,
  Pic(X)) = 2\rho_X.
\]

(ii) \(H^j(X,E)\) is torsion for \(j>0\) and zero for \(j>
  1 +  2~{\rm dim}~X\). 

(iii) \(H^j(X,E)\) is finite for \(j \ge 3\). 

(iv) \(H^1(X,E)\) is finite. 

(v) \(H^2(X,E)\) is co-finite of corank \(2\rho_X\). 

(vi) The \(\hat{\Z}\)-modules \(H^i(X, TE) = \prod_lH^i(X, T_lE)\)
(all primes \(l\) including \(l=p\)) are finite for \(i \neq 1,2\).\medskip 

\noindent (b) The Weil-\'etale cohomology groups \(H^i_W(X,E)\) of \(E\) are
finitely generated. The rank of \(H^i_W(X,E)\) is zero for
\(i\neq0,1\)  and \(2\rho_X\) for \(i=0,1\).
\end{thm} 

\nd {\em Remark.} If the N\'eron-Severi group \(NS(X)\) is torsion-free and
  the Picard scheme \(Pic_X\) is smooth, then \(H^1(X,E)\) is
  isomorphic to the finite group \(\Ext^1_{\F_q} (A,E)\)
  \cite[Thm.~1]{mi2}.\qed

\begin{proof} Part (b) will be proved later (as part of Theorem
  \ref{konege}). We prove part (a).

(i) This is well known (Mordell-Weil theorem);
cf. for instance, \cite{la}.

(ii) It is straightforward that  
 \(H^j(X,E)\) and \(H^j(\bar{X}, E)= H^j_{W}(\bar{X}, E)\)
 are  torsion for \(j>0\) using the fact that \(E\) over \(X\) is the
 Neron model of \(E\) over \(\F_q(X)\);
cf. \cite[p.~41]{breen}.  Consider the Kummer sequence
\begin{equation}\label{kummer}
0 \to~ _nE \to E \xrightarrow{n} E \to 0, \quad (n >0)
\end{equation} 
which is exact in the \'etale (for \(n\) coprime to \(p\)) and flat
topologies. When \((n,p)=1\), the 
identity \(_nE = (\Z/{n\Z})^{2}\) over \(\F\) tells us that \(H^j(\bar{X},~
_nE)\) are finite (zero for \(j>2~{\rm dim} X\)). This implies (ii) up
to \(p\)-torsion  using the  spectral sequence
\begin{equation}\label{ets}
H^r(\Gamma_0, H^s(\bar{X}, E)) \Rightarrow
H^{r+s}(X, E).\end{equation} 
By  \cite[VI, Rmk.~1.5]{mi2}, \(H^j_{fl}({X},~_pE) =0\)  for
\(j>2+{\rm dim}~X\). As \(E\)  is smooth, \(H^j_{fl}({X}, E) =
H^j({X},E)\) \cite[p.~92]{mi2} and (ii) follows.

(iii) It suffices to  consider the non-\(p\) torsion. The
\(p\)-part can be treated similarly, given the sequence
\[
0 \to H^{j-1}({X}, E)\otimes\Z_{p} \to  H^{j}( {X}, T_{p}
E) \to T_{p}H^j( {X}, E) \to 0
\]and the results of the
previous section on \(H^*(X, T_pE)\).  
Though the proof of the non-\(p\)-part is standard
(cf. \cite[Cor. 6.4]{mi4}), nevertheless we recall it for the
convenience of the reader.   

Let \(\gamma\) be the Frobenius automorphism on the Tate module
\(T_{\ell}E\) of \(E\); here \(\ell\) is a prime distinct from
\(p\). By the Riemann hypothesis, \((1-\gamma)\)  is a 
quasi-isomorphism of \(H^j(\bar{X},T_{\ell}E) = H^j(\bar{X},
\Z_{\ell}) \otimes T_{\ell} E\) for \(j\neq1\); thus,
if \( j \neq 1\), then \(H^j(\bar{X}, T_{\ell}E)^{\Gamma}\) and \(H^j(\bar{X},
T_{\ell}E)_{\Gamma}\) are finite. From 
\[
0 \to H^{j-1}(\bar{X}, T_{\ell}E)_{\Gamma} \to H^j(X, T_{\ell}E) \to
H^j(\bar{X}, T_{\ell}E)^{\Gamma} \to 0,
\]we obtain that \( H^j(X, T_{\ell}E) \) is finite for \(j \neq1,2\).  
Now, (\ref{kummer})  provides the exact sequence \begin{equation}\label{etat}
0 \to H^{j-1}({X}, E)\otimes\Z_{\ell} \to  H^{j}( {X}, T_{\ell}
E) \to T_{\ell}H^j( {X}, E) \to 0,\end{equation}whereby, for \(j >2\), \(H^{j}( {X}, T_{\ell}
E)\) is finite and hence \(T_{\ell}H^j( {X}, E)\), being torsion-free,
is zero. As the \(\ell\)-torsion of \(H^j(X,E)\) is finite, the
\(\ell\)-primary subgroup \(H^j(X,E)(\ell)\) itself is finite for \(j
>2\) and isomorphic to \(H^{j+1}(X, T_{\ell}E)\). For the finiteness
of non-\(p\)-part of \(H^j(X,E)\) for \(j\ge 3\), we need the following lemma
whose formulation and proof are as in \cite[Lem.~1.1, pp.176-7]{sl1984}. 
\begin{lem} Write \(P_j(t) = det(1-\gamma t)\) for the characteristic
  polynomial of \(\gamma\) on \(H^j(\bar{X}, \Z_{\ell})\). Similarly,
  we define \(Q_j(t)\) for \(H^j(\bar{X}, T_{\ell}E)\). For
  all \(j\neq 0,1\), we have
\[
[H^j(\bar{X}, E(\ell))^{\Gamma}] = [H^{j+1}(\bar{X},
T_{\ell}E)^{\Gamma}]. |P_j(q^{-\half})|^{-2}_{\ell},
\]where \(|- |_{\ell}\) is the absolute value normalized so that
\(|\ell|_{\ell} = {\ell}^{-1}\). 

One has \(|P_j(q^{-\half}t)|^2_{\ell} = |Q_j(t)|_{\ell}\). \end{lem}

As in \cite[(c), pp.181]{sl1984}, for \(j >2\), \(H^j(X,E)(\ell)\) is
trivial unless \(\ell\) divides \(P_j(q^{-\half})\) or \(H^*(\bar{X},
\Z_{\ell})\) has torsion. By O.~Gabber  \cite{ga},  \(H^j(X,E)(\ell)\) is
nontrivial only for a finite number of primes \(\ell\) thereby showing
the finiteness of the non-\(p\) part of \(H^j(X,E)\) for \(j >2\).  

(iv) Recall that the Hochschild-Serre spectral sequence
\[
H^i(\Gamma, H^j(\bar{X}, T_{\ell}E)) \Rightarrow H^{i+j}(X, T_{\ell}E)
\]yields the
 sequence
\[0 \to H^0(\bar{X}, T_{\ell}E)_{\Gamma} \to H^1({X}, T_{\ell}E) \to
H^1(\bar{X}, T_{\ell}E)^{\Gamma} \to 0.  
\]But
\[
H^1(\bar{X}, T_{\ell}E)^{\Gamma} = \Hom(T_{\ell}A, T_{\ell}E)^{\Gamma}
\xleftarrow{\cong} \Hom_{\F_q} (A,E)\otimes\Z_{\ell}. 
\]As \( H^0(\bar{X}, T_{\ell}E)_{\Gamma}\) is finite, the rank of \(
H^1({X}, T_{\ell}E)\) is \(2\rho_X = {\rm rank}~H^1(\bar{X},
T_{\ell}E)^{\Gamma}\).  The sequence
\[
0 \to H^0(X, E) \otimes\Z_{\ell} \to H^1(X, T_{\ell}E) \to
T_{\ell}H^1(X,E) \to 0,
\]given (i), shows that \(T_{\ell}H^1(X,E)\) has rank zero
and, being torsion-free, is zero. This proves (iv). 

(v) A similar argument as in (iv) using that the rank of
\(H^1(\bar{X}, T_{\ell}E)_{\Gamma}\) is \(2\rho_X\) shows that the
rank of \(T_{\ell}H^2(X,E)\) is \(2\rho_X\) thereby proving (v). 

(vi) was proved in the course of the proof of (iii). 
\end{proof}

We can now formulate an adelic version of Theorem \ref{bures}. 
\begin{thm}\label{adeles} The cohomology groups \(H^i(N^{\bu})\) of the complex
\[
N^{\bu}:\quad \cdots H^i(X, TE) \xrightarrow{\alpha_i} H^{i+1}(X,TE)
\to \cdots
\]are finite; the maps \(\alpha_i\)  are induced by cup-product with
the generator of \(H^1({\rm Spec}~\F_q, \hat{\Z})\). One has  
\[
c_X^2 = q^{\chi(X, \mathcal O_X)} \prod  ([H^i(N^{\bu})])^{(-1)^i}. 
\]
\end{thm} 

\begin{proof} The first statement needs verification only in the cases
  \(i=1,2\), given (vi) of Theorem  \ref{manins}. We need to check
  that the kernel and cokernel of \(H^1(X,TE) \to H^2(X,TE)\) are
  finite. The \(p\)-part is proved in Theorem \ref{bures}. The
  non-\(p\) part is straightforward, given the previous lemma and the obvious
  variant of the non-\(p\)-part of \cite[Lemma~6.2]{mi4} (or the
  \'etale case of \cite[Lemma~4]{mi1}); see the
  proof of \cite[Theorem~0.1]{mi4}.\end{proof}

\subsection{Proof of Theorem \ref{konege}}  
Recall \cite[\S8]{sl2003} the generator \(\theta\) of \(H^1_W({\rm Spec}~\F_q,
\Z) ={\Hom}(\Z, \Z) =\Z\). 
It is clear that \(H^0_W(X,E) = H^0(X,E)\)  is
finitely generated.  

The groups \(H^i(\bar{X}, E)\) are torsion for \(i>0\)
\cite[p.~41]{breen}.  By Theorem \ref{manins} and
\cite[Prop.~2.3]{sl2003}, \(H^i(X,E) \cong 
H^i_W(X,E)\)  for \(i\ge 3\) and thus proves (ii). To prove (i), it
remains to show that \(H^1_W(X,E)\) and \(H^2_W(X,E)\) are finitely
generated. 
Now, \cite[Prop.~2.3]{sl2003} provides a
spectral sequence
\[H^r(\Gamma_0, H^s_W(\bar{X}, E)) \Rightarrow
H^{r+s}_W(X, E)
\]with a map from (\ref{ets}) of spectral sequences. This gives exact sequences
\[
0 \to H^2(\Gamma, H^0(\bar{X}, E)) \to H^2(X,E) \to H^2_W(X,E) \to 0,
\]and
\[
0 \to  H^1(X,E) \to H^1_W(X,E) \to \frac{H^1(\Gamma_0,
  H^0(\bar{X},E))}{H^1(\Gamma, H^0(\bar{X},E))}\to 0.
\]
As \(H^0(\bar{X},E)\) differs from \(Y:= \Hom_{\F}(\bar{A},\bar{E})\)
by torsion, we have 
\[
H^2(\Gamma, H^0(\bar{X}, E)) \cong H^1(\Gamma_0, H^0(\bar{X},E))
\otimes{\frac{\Q}{\Z}} \cong H^1(\Gamma_0, Y)\otimes{\frac{\Q}{\Z}},
\]using \cite[Lemma~1.2]{sl2003}. The sequence of \(\Gamma_0\)-modules
\(0\to Y \xrightarrow{n} Y \to Y/{nY} \to 0\) gives
\[
H^1(\Gamma_0,Y)\otimes{\frac{\Q}{\Z}} 
\cong H^1(\Gamma_0, Y\otimes{\frac{\Q}{\Z}}).
\]and, by
Tate's theorem \cite{Tatendo}, the corank of \(H^1(\Gamma_0, Y\otimes{\frac{\Q}{\Z}})\) is \(2\rho_X\).  This proves that
\(H^2_W(X,E)\) is finite, the rank of \(H^1_W(X,E)\) is \(2\rho_X\)
thereby finishing the proof of (i), (iii) and (iv). 

  Part (v) follows from Theorem \ref{adeles} and the isomorphisms
  \(H^i_W(X,E)\otimes\hat{\Z} 
\cong H^{i+1}(X, TE)\) (see (\ref{etat}); the nontrivial cases are \(i=0,1\)) which are compatible with
the maps  \(\cup \theta\) 
and \(\alpha_i\). In more detail, combining \(TH^1(X,E)=0\) 
(Theorem \ref{manins} (iv)), the isomorphism \(H^0(X,E) =
H^0_W(X,E)\), and (\ref{etat}) proves the case \(i=0\). The exactness
of (\ref{kummer}) in the Weil-\'etale or Weil-flat
topology gives an exact sequence
\[
0 \to H^i_W(X,E)\otimes{\hat{\Z}} \to H^{i+1}_W(X,TE) \to
TH^{i+1}_W(X,E) \to 0,
\]where the first map is an isomorphism because \(H^*_W(X,E)\) are
finitely generated abelian groups and so \(TH^*_W(X,E)=0\). As
\(H^*(X,TE) \xrightarrow{\sim}H^*_W(X,TE)\), this finishes the
proof of Theorem \ref{konege}.\qed

\subsection*{Acknowledgements} I heartily thank P. Cartier,
C. Deninger, O. Gabber, 
M. Kontsevich, 
S. Lichtenbaum, Yu. I. Manin, and J. Milne for comments and
suggestions. I am grateful to the referees for many constructive
comments. This work was begun at Max-Planck
Institut f\"{u}r Mathematik, Bonn and completed at Institut des
Hautes \'{E}tudes Scientifiques; I thank both institutions for their
hospitality and 
support.

\end{document}